\documentclass[graybox]{svmult}

\usepackage{mathptmx}       
\usepackage{helvet}         
\usepackage{courier}        
\usepackage{type1cm}        

\usepackage{graphicx}        
\usepackage[bottom]{footmisc}

\usepackage{amssymb}
\usepackage{xypic}

\usepackage{enumerate}

\DeclareMathAlphabet{\mathcal}{OMS}{cmsy}{m}{n}

\newcounter{number}[section]

\newenvironment{nummer}{\refstepcounter{number}{\bf\noindent\arabic{section}.\arabic{number} }}{}

\newcommand{\bn}{\noindent \begin{nummer} \rm}
\newcommand{\en}{\end{nummer}\bigskip}

\newenvironment{ntheorem}{\noindent {\bf Theorem} \it}{}
\newenvironment{nlemma}{\noindent {\bf Lemma} \it}{}
\newenvironment{nprop}{\noindent {\bf Proposition} \it}{}
\newenvironment{ndefn}{\noindent {\bf Definition} }{}
\newenvironment{ncor}{\noindent {\bf Corollary} \it}{}

\newenvironment{nremark}{\noindent {\bf Remark}}{}
\newenvironment{nremarks}{\noindent {\bf Remarks}}{}
\newenvironment{nexamples}{\noindent {\bf Examples} }{}

\newenvironment{nquestion}{\noindent {\bf Question} }{}
\newenvironment{nquestions}{\noindent {\bf Questions} }{}
\newenvironment{nproof}{\noindent {\bf Proof}}{\mbox{}\hfill 
\rule[-.2ex]{.25em}{1.8ex}}
\newenvironment{idea}{\noindent {\bf Idea of proof}}{\mbox{}\hfill 
\rule[-.2ex]{.25em}{1.8ex}}

\begin{document}

\title*{QDQ vs.\ UCT}
\author{Wilhelm Winter}
\institute{Wilhelm Winter \at Mathematisches Institut der Westf\"alischen Wilhelms-Universit\"at M\"unster, Einsteinstra{\ss}e 62, 48149 M\"unster, Deutschland; \email{wwinter@uni-muenster.de}}

%
%
\maketitle

\abstract{
This is a survey of recent progress in the structure and classification theory of nuclear $\mathrm{C}^{*}$-algebras. In particular, I outline how the Universal Coefficient Theorem ensures a positive answer to the quasidiagonality question in the presence of faithful traces. This has strong consequences for the regularity conjecture and the classification problem for separable, simple, nuclear $\mathrm{C}^{*}$-algebras. Moreover, it entails a positive solution to Rosenberg's conjecture on quasidiagonality of reduced $\mathrm{C}^{*}$-algebras of discrete amenable groups. This note is largely based on a joint paper with Aaron Tikuisis and Stuart White.}

\section*{Introduction}

\noindent
Quasidiagonality was defined by Halmos as an external finite dimensional approximation property for sets of operators on a Hilbert space \cite{Hal:BAMS}. Voiculescu studied the notion as a property of (represented) $\mathrm{C}^{*}$-algebras; cf.\ \cite{Voi:Duke,V:IEOT}. 

For nuclear $\mathrm{C}^{*}$-algebras quasidiagonality can be expressed as an embedding property via the Choi--Effros lifting theorem \cite{CE:Ann}. Another connection to amenability was disclosed by Hadwin in \cite{Hadwin:JOT} and---in an appendix to the same paper---by Rosenberg who observed that discrete group $\mathrm{C}^{*}$-algebras can only be quasidiagonal for amenable groups. The converse statement became known as Rosenberg's conjecture. 

Elliott's programme to classify nuclear $\mathrm{C}^{*}$-algebras by $\mathrm{K}$-theory data first featured quasidiagonality in the proof of Kirchberg's embedding theorem \cite{Kir:ICM} (where it entered through Voiculescu's homotopy invariance theorem \cite{Voi:Duke}), which then led up to the celebrated Kirchberg--Phillips classification of simple, purely infinite, nuclear $\mathrm{C}^{*}$-algebras; see \cite{KP:Crelle, R:Book}. 

In the stably finite case, the relevance of quasidiagonality was first marked by Popa's local quantisation of \cite{Popa:PJM}, which then inspired Lin to define TAF algebras, thus leading up to the whole TAF/TAI classification machinery (which recently culminated in the spectacular \cite{GLN:arXiv} and \cite{EGLN:arXiv} by Elliott, Gong, Lin, and Niu). 

The quasidiagonality question (QDQ) of Blackadar and Kirchberg asks whether all separable, stably finite, nuclear $\mathrm{C}^{*}$-algebras are quasidiagonal; see \cite{BlaKir:MathAnn}. In the 2000's, the interest in the interplay between quasidiagonality and nuclearity was renewed in Brown's \cite{B:MAMS} (which takes a measure theoretic or tracial point of view) and in \cite{KirWin:dr} (which establishes a further link to topological approximation properties via the decomposition rank, a noncommutative covering dimension). The latter notion was generalised to nuclear dimension in \cite{WinZac:dimnuc}, and it became an intriguing problem to characterise the difference between decomposition rank and nuclear dimension. 

A massive breakthrough was achieved by Matui and Sato in \cite{MS:DMJ}: In the simple, monotracial and $\mathcal{Z}$-absorbing case (where $\mathcal{Z}$ denotes the Jiang--Su algebra from \cite{JS:AJM}; see also \cite{RorWin:Z-revisited}) quasidiagonality implies finite decomposition rank; in the absence of quasidiagonality, by  \cite{SWW:Invent} one always has finite nuclear dimension. For larger trace spaces, one has to consider Brown's more refined notion of quasidiagonality of traces; cf.\ \cite{BBSTWW:arXiv}. 

The latter  condition also appears in \cite{EGLN:arXiv} which amounts to the classification of all separable, unital, simple $\mathrm{C}^{*}$-algebras with finite nuclear dimension, which satisfy the Universal Coefficient Theorem from \cite{RS:DMJ} (UCT for short; satisfying it is equivalent to being $\mathrm{KK}$-equivalent to a commutative $\mathrm{C}^{*}$-algebra) and for which all traces are quasidiagonal. Although of stunning generality, the last two hypotheses in this result do remain mysterious. The UCT holds for $\mathrm{C}^{*}$-algebras of amenable groupoids by \cite{Tu:KT} and is also ensured by a bootstrap type condition; therefore in concrete situations of interest it can usually be confirmed (see, for example, \cite{EckGil:irr}). To date it is unknown whether \emph{all} separable nuclear $\mathrm{C}^{*}$-algebras satisfy the UCT. Quasidiagonality---due to its poor permanence properties---is usually also hard to verify even in very concrete and geometric setups; see for example \cite{Pim} or \cite{Lin:qd}. 

Remarkably, neither of these two problems occur in the von Neumann algebra context of Connes' celebrated classification of injective factors (which remains an inspiration for the classification of nuclear $\mathrm{C}^{*}$-algebras---at a philosophical, but also at a technical level, since the various available proofs in \cite{C:Ann, H:JFA, Popa:JOT} have revealed important insights for $\mathrm{C}^{*}$-algebraists). 

The situation on the $\mathrm{C}^{*}$-algebra side changed considerably with \cite{TWW}, in which the two problems were linked to the effect that the UCT indeed ensures quasidiagonality under suitable circumstances, with the crucial extra tool  being stable uniqueness theorems as introduced in  \cite{DE:PLMS, Lin:stableJOT}. The result verifies Rosenberg's conjecture and has strong consequences for the structure and classification of simple nuclear $\mathrm{C}^{*}$-algebras, including the regularity conjecture of Toms and myself (see \cite{ET:BAMS, W:Invent1}). In particular, the classification of separable, unital, simple $\mathrm{C}^{*}$-algebras which satisfy the UCT and have finite nuclear dimension is now complete. Moreover, the nuclear dimension hypothesis may be replaced by $\mathcal{Z}$-absorption or by strict comparison of positive elements, at least under additional conditions on the space of tracial states.

In this note I give a survey of the result, its proof and its corollaries. I try to motivate it and put it into a larger context---especially in view of the remaining open problems in the area---by specialising to the case of strongly self-absorbing $\mathrm{C}^{*}$-algebras.

In Section~\ref{preliminaries-qd} I recall the concept of and some elementary facts about quasidiagonality. In Section~\ref{su-UCT-preliminaries} I state the stable uniqueness theorems and outline the role of the UCT for the controlled version. Section~\ref{QDQ-section} looks at strongly self-absorbing $\mathrm{C}^{*}$-algebras and approaches various versions of the quasidiagonality question from this point of view. Section~\ref{main-section} contains the main result and its consequences for structure and classification of nuclear $\mathrm{C}^{*}$-algebras and for Rosenberg's conjecture. The proof of the main theorem is outlined in Section~\ref{outline-section}. Finally, Section~\ref{open-section} collects a number of open problems, mostly on quasidiagonality and on the UCT. 

Most of the contents of the paper are based on \cite{TWW}, except for  those of Sections~\ref{QDQ-section} and \ref{open-section} which have not previously been published in this form. I would like to thank Ilijas Farah, David Kerr, Aaron Tikuisis, Stuart White and the referee for some very helpful comments on an earlier version. 

This note was written during a research stay at the Mittag-Leffler Institute to which I am indebted for their hospitality. I am grateful to the Deutsche Forschungsgemeinschaft for their support through the collaborative research centre SFB 878. And I would like to thank the organisers of the 2015 Abel Symposium for their kind invitation to this marvellous conference.

\section{Quasidiagonality}
\label{preliminaries-qd}

\bn
\label{qd-def}
Quasidiagonality was originally defined by Halmos for sets of operators on a Hilbert space; an abstract $\mathrm{C}^{*}$-algebra is called quasidiagonal if it has a faithful representation the image of which is quasidiagonal in this sense. For our purposes, it will be useful to take Voiculescu's characterisation from \cite{Voi:Duke} as a definition:

\bigskip

\begin{ndefn}
A $\mathrm{C}^{*}$-algebra $A$ is quasidiagonal if, for every finite subset $\mathcal{F}$ of $A$ and $\varepsilon > 0$, there exist $k \in \mathbb{N}$ and a completely positive contractive (c.p.c.) map $\psi:A \to M_{k}$ such that 
\[
\|\psi(ab) - \psi(a) \psi(b)\|,\;  \|a\| - \|\psi(a)\| < \varepsilon
\]
for all $a,b \in \mathcal{F}$.
\end{ndefn}
\en

\bn
\label{qd-prop}
With the aid of the Choi--Effros lifting theorem it is not hard to see that for separable nuclear $\mathrm{C}^{*}$-algebras quasidiagonality can be rephrased in terms of embeddings into $\mathcal{Q}_{\omega}$, where $\mathcal{Q}$ is the universal UHF algebra, $\omega \in \beta \mathbb{N} \setminus \mathbb{N}$ is a free ultrafilter, and $\mathcal{Q}_{\omega} = \prod_{\mathbb{N}} \mathcal{Q}/ \{(x_{n})_{\mathbb{N}} \mid \lim_{n \to \omega} x_{n} = 0\}$. 

\bigskip

\begin{nprop}
A separable nuclear $\mathrm{C}^{*}$-algebra $A$ is quasidiagonal if and only if there is an embedding $A \hookrightarrow \mathcal{Q}_{\omega}$.
\end{nprop}
\en

\bn
\label{unital-qd}
\begin{nremark}
Note that in this characterisation unitality is not an issue: $A$ as above is quasidiagonal if and only if the smallest unitisation $A^{\sim}$ is quasidiagonal if and only if there is a unital embedding $A^{\sim} \hookrightarrow \mathcal{Q}_{\omega}$. (This last statement uses the fact that $\mathcal{Q}$ is self-similar in the sense that every unital hereditary $\mathrm{C}^{*}$-subalgebra of $\mathcal{Q}$ is isomorphic to $\mathcal{Q}$.)
\end{nremark}
\en

\bn
Brown considered in \cite{B:MAMS} the refined notion of \emph{quasidiagonal traces}; this is based on Voiculescu's observation \cite{V:IEOT} that a unital, separable, quasidiagonal $\mathrm{C}^{*}$-algebra always has at least one tracial state which is picked up by a sequence of (unital) quasidiagonal approximations as in Definition~\ref{qd-def}. 

\bigskip

\begin{ndefn}
A tracial state $\tau_{A} \in \mathrm{T}(A)$ on a $\mathrm{C}^{*}$-algebra $A$ is quasidiagonal if, for every finite subset $\mathcal{F}$ of $A$ and $\varepsilon > 0$, there exist $k \in \mathbb{N}$ and a c.p.c.\ map $\psi:A \to M_{k}$ such that 
\[
\|\psi(ab) - \psi(a) \psi(b)\|,\;  | \tau_{A}(a) - \tau_{M_{k}} \circ \psi(a) | < \varepsilon
\]
for all $a,b \in \mathcal{F}$.
\end{ndefn}
\en

\bn
\label{qd-traces-prop}
Similar to Proposition~\ref{qd-prop}, one can characterise quasidiagonality of traces in terms of maps to $\mathcal{Q}_{\omega}$.

\bigskip

\begin{nprop}
Let $A$ be a separable nuclear $\mathrm{C}^{*}$-algebra. A tracial state $\tau_{A} \in \mathrm{T}(A)$ is quasidiagonal if and only if there is a $^{*}$-homomorphism $\pi:A \to \mathcal{Q}_{\omega}$ such that $\tau_{A} = \tau_{\mathcal{Q}_{\omega}} \circ \pi$, where $\tau_{\mathcal{Q}_{\omega}}$ is the unique tracial state on $\mathcal{Q}_{\omega}$ given by $\tau_{\mathcal{Q}_{\omega}}([(x_{n})_{\mathbb{N}}]) = \lim_{n \to \omega} \tau_{\mathcal{Q}}(x_{n})$. 
\end{nprop}
\en

\section{Stable uniqueness and the Universal Coefficient Theorem}

\label{su-UCT-preliminaries}

\bn
\label{su-original}
It is a common problem in classification to decide when two morphisms are (approximately unitarily) equivalent, provided there is no obvious $\mathrm{K}$-theoretic obstruction. Stable uniqueness theorems provide partial solutions by establishing `local almost' unitary equivalences after adding the same map on both sides with sufficiently large multiplicity. (Very roughly speaking, this could perhaps be interpreted as a fine-tuned version of Voiculescu's theorem \cite{Voi:WeylVN}.) The problem has been studied extensively by Lin (see \cite{Lin:stableJOT}, for example); for our purposes a version of Dadarlat and Eilers---which we state here in a simplified form for brevity---will be particularly useful.

\bigskip

\begin{ntheorem}
{\rm \cite[Theorem~4.5]{DE:PLMS} }
Let $C,B$ be unital $\mathrm{C}^{*}$-algebras with $C$ separable and nuclear. Let $\iota:C \to B$ be a unital $^{*}$-homomorphism which is totally full, i.e., for every nonzero $c \in C$ the element $\iota(c)$ generates $B$ as an ideal. Let $\varphi,\psi: C \to B$ be unital $^{*}$-homomorphisms with the same induced $\mathrm{KK}$-class. 

Then, for every finite subset $\mathcal{G} \subset C$ and $\delta>0$ there are $n \in \mathbb{N}$ and a unitary $u \in M_{n+1}(B)$ such that
\[
\| u (\varphi(c) \oplus \iota^{\oplus n}(c)) u^{*} - (\psi(c) \oplus \iota^{\oplus n}(c)) \| < \delta, \quad c \in \mathcal{G}.
\]  
\end{ntheorem}
\en

\bn
\begin{nremark}
For us it will be important that the domain algebra $C$ in the theorem above need not be simple. In fact, we will apply it to the unitisation of the suspension of the algebra we are actually interested in. For the target we mostly care about the algebra $\mathcal{Q}_{\omega}$; for technical reasons we also consider algebras like $\prod_{\mathbb{N}} \mathcal{Q}_{\omega}$.
\end{nremark}
\en

\bn
For the proof of our main result, we need a refined version of Theorem~\ref{su-original}, since for our method it is vital that $n$ can be chosen independent of the specific maps $\varphi, \psi$ and especially $\iota$ (it may of course still depend on $C$, $\mathcal{G}$ and $\delta$). Such a result was already provided in \cite[Theorem~4.12]{DE:PLMS}, at least for simple domains $C$. However, the proof carries over to the nonsimple situation as long as one retains some control over the fullness of $\iota$. This is done in \cite{TWW} in terms of a control function $\Delta: C_{+}^{1} \setminus \{0\} \to \mathbb{N}$ (where $C_{+}^{1}$ is the unit ball of the positive elements of $C$ and $\Delta(c)$ is the smallest $k \in \mathbb{N}$ such that there are contractions $x_{1}, \ldots, x_{k} \in B$ with $1_{B} = \sum_{i=1}^{k} x_{i}^{*} \iota(c) x_{i}$ ). For a given $^{*}$-homomorphism $\iota: C \to B$ one can define an (a priori possibly infinite) function $\Delta$ in terms of (inverses of) tracial values of $\iota(c)$ for $c \in C_{+}^{1} \setminus \{0\}$; provided $B$ has strict comparison of positive elements and the control function remains finite, the map $\iota$ will then be $\Delta$-full---see \cite{TWW} for details.

The proof of \cite[Theorem~4.12]{DE:PLMS} works by assuming the statement is false and then producing a sequence of pairs of counterexamples (i.e., pairs of maps as in Theorem~\ref{su-original}, but with no uniform bound on $n$) for contradiction. The crux is it to assemble the two resulting sequences of maps into two single maps (with larger target algebras) and at the same time keeping control over their $\mathrm{KK}$-classes. In our situation, the individual maps will be zero-homotopic, and the problem is it to decide when maps of the form
\begin{equation}
\label{KK-product-injective}
\textstyle
\mathrm{KK}\big(C,\prod_{\mathbb{N}}B_{n}\big) \to \prod_{\mathbb{N}} \mathrm{KK}(C,B_{n})
\end{equation}
are injective. The issue is that a sequence of homotopies doesn't necessarily give rise to a (continuous) homotopy of sequences, since the parameter speeds might increase quickly. It is for this purpose---and for this purpose only---that the Universal Coefficient Theorem enters the stage.    
\en

\bn
\begin{ndefn}{\rm \cite[Theorem~1.17]{RS:DMJ}, \cite{B:KThy} }
A separable $\mathrm{C}^{*}$-algebra $A$ is said to satisfy the Universal Coefficient Theorem (UCT) if the sequence
\[
0 \to \mathrm{Ext}(\mathrm{K}_{*}(A),\mathrm{K}_{*+1}(B)) \to \mathrm{KK}(A,B) \to \mathrm{Hom}(\mathrm{K}_{*}(A),\mathrm{K}_{*}(B)) \to 0
\]
is exact for every $\sigma$-unital $\mathrm{C}^{*}$-algebra $B$.
\end{ndefn}

\bigskip

Both maps in the sequence above can be made explicit, using the natural identification $\mathrm{KK}(A,B) \cong \mathrm{Ext}^{-1}(A, \mathcal{C}_{0}(\mathbb{R}) \otimes B \otimes \mathcal{K})$: Given such an extension, the right hand map collects the boundary maps of the associated six-term exact sequence in $\mathrm{K}$-theory. For trivial boundary maps, the six-term exact sequence splits into two extensions of abelian groups; the UCT requires the left-hand map to be the inverse of this assignment. 
\en

\bn
It is known that the separable $\mathrm{C}^{*}$-algebras satisfying the UCT are precisely the ones which are $\mathrm{KK}$-equivalent to abelian $\mathrm{C}^{*}$-algebras; see \cite{RS:DMJ}, \cite[Proposition~5.3]{Skandalis:KT} (or \cite[Theorem~23.10.5]{B:KThy}). The closure properties of the class of \emph{nuclear} UCT $\mathrm{C}^{*}$-algebras are so strong that to date nobody has managed to find a nuclear example outside it. The UCT problem reads as follows.

\bigskip

\begin{nquestion}
Does every separable nuclear $\mathrm{C}^{*}$-algebra satisfy the UCT?
\end{nquestion}  
\en

\bn
\label{su-controlled}
In our situation, the UCT provides just enough information to make the map of (\ref{KK-product-injective}) injective; see \cite{DE:PLMS} and \cite{TWW} for details. Let us state the resulting `controlled' stable uniqueness theorem, again for brevity in a simplified version. (See \cite{TWW} for the complete statement; this involves $\mathrm{K}$-theory with coefficients, which here boils down to just ordinary $\mathrm{K}$-theory by $\mathcal{Q}$-stability.)

\bigskip

\begin{ntheorem}
Let $C$ be a separable, unital, nuclear $\mathrm{C}^{*}$-algebra satisfying the UCT. Let $\Delta: C^{1}_{+} \setminus \{0\} \to \mathbb{N}$ be a control function, let $\mathcal{G} \subset C$ be a finite subset and let $\delta >0$. Then there exists $n \in \mathbb{N}$ such that for any unital $\Delta$-full $^{*}$-homomorphism $\iota: C \to B$ (where $B$ is of the form $\mathcal{Q}$, $\mathcal{Q}_{\omega}$ or $\prod_{\mathbb{N}} \mathcal{Q}_{\omega}$), and any unital $^{*}$-homomorphisms $\varphi,\psi: C \to B$ with $\mathrm{K}_{*}(\varphi) = \mathrm{K}_{*}(\psi)$, there is a unitary $u \in M_{n+1}(B)$ such that
\[
\| u (\varphi(c) \oplus \iota^{\oplus n}(c)) u^{*} - (\psi(c) \oplus \iota^{\oplus n}(c)) \| < \delta, \quad c \in \mathcal{G}.
\] 
\end{ntheorem} 
\en

\section{QDQ: a strongly self-absorbing point of view}
\label{QDQ-section}

\bn
The quasidiagonality question, also known as Blackadar--Kirchberg problem, was posed in \cite{BlaKir:MathAnn}:

\bigskip
\noindent
{$\bf QDQ$} \enspace {\it Is every separable, stably finite, nuclear $\mathrm{C}^{*}$-algebra quasidiagonal?}

\bigskip
\noindent
Note that by Remark~\ref{unital-qd}, {$\bf QDQ$} is equivalent to {${\bf QDQ}_{1}$}, the respective question for unital $\mathrm{C}^{*}$-algebras. 

Likewise, one could add simplicity or both unitality and simplicity to the hypotheses; we denote the resulting questions by {${\bf QDQ}_{\mathrm{simple}}$} and {${\bf QDQ}_{\mathrm{simple},1}$}. At this point I am not aware of any (but the obvious) formal implications between {$\bf QDQ$}, {${\bf QDQ}_{\mathrm{simple}}$} and {${\bf QDQ}_{\mathrm{simple},1}$}, although it doesn't seem unlikely that the general case can be reduced to the simple situation (for example via some Bernoulli type crossed product construction as in \cite{ORS:GAFA}?).
\en

\bn
Although the quasidiagonality question has been around for a long time, until recently its role for the structure and classification theory of nuclear $\mathrm{C}^{*}$-algebras has remained somewhat obscure, and its nature and complexity is still hard to gauge. Maybe most conspiciously, it seems to be a very finite problem, in stark contrast to other structural phenomena, which often exhibit a dichotomy between finite and infinite situations, such as the interplay between Lin's TAF classification and Kirchberg--Phillips classification, the roles of the Jiang--Su algebra $\mathcal{Z}$ and the Cuntz algebra $\mathcal{O}_{\infty}$, and Toms' and R{\o}rdam's topologically high-dimensional examples \cite{T:Ann, R:Acta}. 
\en

\bn
In order to shed some light on the quality of the problem, let us be more restrictive and formulate the quasidiagonality question for strongly self-absorbing $\mathrm{C}^{*}$-algebras. These have been abstractly defined in \cite{TW:TAMS}, but the concept goes back to Effros--Rosenberg \cite{EffRos:PJM} in the $\mathrm{C}^{*}$-algebra setting, which in turn was inspired by the von Neumann algebra context where it was crucial for McDuff's \cite{M:PLMS} and for Connes' celebrated classification of injective factors \cite{C:Ann}.
\en

\bn
\begin{ndefn}
A separable unital $\mathrm{C}^{*}$-algebra $\mathcal{D} \neq \mathbb{C}$ is strongly self-absorbing, if there is an isomorphism $\varphi: \mathcal{D} \to \mathcal{D} \otimes \mathcal{D}$ which is approximately unitarily equivalent to the first factor embedding, i.e., there is a sequence of unitaries $(u_{n})_{\mathbb{N}}$ in $\mathcal{D} \otimes \mathcal{D}$ such that $u_{n} (d \otimes 1_{\mathcal{D}})u_{n}^{*} \stackrel{n \to \infty}{\longrightarrow} \varphi(d)$ for all $d \in \mathcal{D}$. 
\end{ndefn}
\en

\bn
It was shown in \cite{EffRos:PJM} that strongly self-absorbing $\mathrm{C}^{*}$-algebras are always simple and nuclear; from results of Kirchberg it follows that they are either purely infinite or stably finite and in this case there is a unique tracial state (cf.\ \cite{TW:TAMS}). By \cite{W:JNCG} they are always $\mathcal{Z}$-stable, i.e.\ they absorb $\mathcal{Z}$ tensorially. 
\en

\bn
For strongly self-absorbing $\mathrm{C}^{*}$-algebras tensorial absorption can be characterised in terms of (exact or approximate) unital embeddings: A separable unital $\mathrm{C}^{*}$-algebra $A$ absorbs the strongly self-absorbing $\mathrm{C}^{*}$-algebra $\mathcal{D}$ precisely if $\mathcal{D}$ embeds unitally into the commutant of $A$ inside its ultrapower, $A_{\omega} \cap A'$ (the criterion can also be phrased for nonunital $A$). The proof is based on an Elliott intertwining argument; see \cite[Theorem~7.2.2]{R:Book} or \cite[Theorem~2.3]{TW:TAMS} (these are stated for sequence algebras instead of ultrapowers, but the proofs are essentially the same). For a strongly self-absorbing target $\mathcal{E}$, a separable subspace of the ultrapower $\mathcal{E}_{\omega}$ can be unitarily conjugated to a subspace of the relative commutant $\mathcal{E}_{\omega} \cap \mathcal{E}'$. As a consequence, two strongly self-absorbing $\mathrm{C}^{*}$-algebras $\mathcal{D}$ and $\mathcal{E}$ are isomorphic if and only if there are unital embeddings $\mathcal{D} \hookrightarrow \mathcal{E}$ and $\mathcal{E} \hookrightarrow \mathcal{D}$ if and only if there are unital embeddings $\mathcal{D} \hookrightarrow \mathcal{E}_{\omega}$ and $\mathcal{E} \hookrightarrow \mathcal{D}_{\omega}$; see \cite{EffRos:PJM} and \cite{TW:TAMS}. Upon combining this with Proposition~\ref{qd-prop}, we see that a (necessarily finite) strongly self-absorbing $\mathcal{D}$ is quasidiagonal if and only if $\mathcal{D} \otimes \mathcal{Q} \cong \mathcal{Q}$. Moreover, Kirchberg's embedding theorem yields that $\mathcal{O}_{2}$ absorbs any other strongly self-absorbing $\mathrm{C}^{*}$-algebra, $\mathcal{O}_{2} \otimes \mathcal{D} \cong \mathcal{O}_{2}$. (This last statement holds in much greater generality.)
\en

\bn
\label{ssa-examples}
\begin{nexamples}
The chart below contains all known strongly self-absorbing $\mathrm{C}^{*}$-algebras. Here, $\mathrm{UHF}^{\infty}$ stands collectively for UHF algebras of infinite type. (The universal UHF algebra $\mathcal{Q}$ is one of them; we list it separately to emphasise its role as a `semifinal' object.) An arrow means `embeds unitally into' or equivalently `is tensorially absorbed by'. 

\vspace{0.5em}

\begin{center}
\begin{tabular}{ccc}
&& $\mathcal{O}_{2}$\\
\\
& $\nearrow$& $\uparrow$ \\
\\
$\mathcal{Q}$ &$ \rightarrow $& $\mathcal{Q} \otimes \mathcal{O}_{\infty}$ \\
\\
$\uparrow$ && $\uparrow$ \\
\\
$\mathrm{UHF}^{\infty}$ & $ \rightarrow $ & $\mathrm{UHF}^{\infty} \otimes \mathcal{O}_{\infty}$ \\
\\
$\uparrow$ && $\uparrow$ \\
\\
$\mathcal{Z}$ & $ \rightarrow $ & $\mathcal{O}_{\infty}$
\end{tabular}
\end{center}
\vspace{0.5em}

\noindent
Arguably the most important question about strongly self-absorbing $\mathrm{C}^{*}$-algebras is whether or not the list above is complete. This makes direct contact with fundamental open problems such as the classification problem, the Toms--Winter conjecture, the UCT problem, or the quasidiagonality question. Even though being strongly self-absorbing is a very restrictive condition, at this point there is no evidence these questions will be substantially easier to answer when restricted to the strongly self-absorbing situation. On the other hand, such a restriction can often bare the problem of merely technical additional complications, and in this way sometimes  disclose its true nature. Occasionally, a solution in the strongly self-absorbing case will then even give us a clue of how to deal with the general situation. This has happened for example in the run-up to \cite{SWW:Invent} and to \cite{TWW}; it is one of the reasons why I like to think of strongly self-absorbing $\mathrm{C}^{*}$-algebras as a microcosm within the world of all nuclear $\mathrm{C}^{*}$-algebras. 
\end{nexamples}
\en

\bn
\label{ssa-initial-final}
It is a crucial feature of the point of view above that questions on the existence or non-existence of examples with certain properties can be rephrased in terms of abstract characterisations of the known examples. For instance, the Jiang--Su algebra $\mathcal{Z}$ was characterised in \cite{W:JNCG} as the uniquely determined initial object in the category of all strongly self-absorbing $\mathrm{C}^{*}$-algebras. (An object in a category is initial, if there is a morphism to every other object. Very often this morphism is also required to be unique; in our situation, this will be the case when using approximate unitary equivalence classes instead of just unital $^{*}$-homomorphisms.) At the opposite end, $\mathcal{O}_{2}$ is the unique final object (i.e., there is a morphism from every other object to $\mathcal{O}_{2}$; as above, this will be unique when using as morphisms approximate unitary equivalence classes of unital $^{*}$-homomorphisms) by Kirchberg's embedding theorem. These are, as Kirchberg once put it, \emph{sociological} characterisations, based on interactions with peer objects. In \cite{DadWin:MS}, it was observed that $\mathcal{O}_{2}$ can also be characterised intrinsically---or \emph{genetically}---as the unique strongly self-absorbing $\mathrm{C}^{*}$-algebra with trivial $\mathrm{K}_{0}$-group. Conspicuously, this characterisation of the \emph{final} object does not require the UCT; in contrast, Kirchberg has shown that the UCT problem is in fact equivalent to the question whether a unital Kirchberg algebra with trivial $\mathrm{K}$-theory is isomorphic to $\mathcal{O}_{2}$ \cite{K:Abel}, and Dadarlat has a parallel result featuring $\mathcal{Q}$ \cite{Dad:UCTremarks}. It is tempting to think of $\mathcal{Q}$ and $\mathcal{Q} \otimes \mathcal{O}_{\infty}$ in a similar way, as `semifinal' objects: $\mathcal{Q}$ is final in the category of all \emph{known} finite strongly self-absorbing $\mathrm{C}^{*}$-algebras, and, more abstractly, also in the category of all quasidiagonal strongly self-absorbing $\mathrm{C}^{*}$-algebras (cf.\ \cite{EffRos:PJM}). One can now turn the tables and interpret this fact as a characterisation of quasidiagonality for strongly self-absorbing $\mathrm{C}^{*}$-algebras in terms of its final object. Similarly, $\mathcal{Q} \otimes \mathcal{O}_{\infty}$ is the final object in the category of all known strongly self-absorbing $\mathrm{C}^{*}$-algebras which are not $\mathcal{O}_{2}$.  Turning tables again one can look at the category of all strongly self-absorbing $\mathrm{C}^{*}$-algebras which embed unitally into $\mathcal{Q} \otimes \mathcal{O}_{\infty}$ and interpret this as a notion of quasidiagonality which also makes sense in the infinite setting, at least in the strongly self-absorbing context.     
\en

\bn
\label{finite-ssa-q}
The strongly self-absorbing version of the quasidiagonality question reads: {\it Is every finite strongly self-absorbing $\mathrm{C}^{*}$-algebra quasidiagonal?} In view of the preceding discussion, we obtain an equivalent formulation as follows:

\bigskip
\noindent 
{${\bf QDQ}_{\mathrm{finite \; s.s.a.}}$}  
{\it Is $\mathcal{D} \otimes \mathcal{Q} \cong \mathcal{Q}$ for every finite strongly self-absorbing $\mathrm{C}^{*}$-algebra $\mathcal{D}$?}

\bigskip
\noindent
Note that this asks whether $\mathcal{Q}$ can be characterised abstractly as the final object in the category of finite strongly self-absorbing $\mathrm{C}^{*}$-algebras. In the above one could specialise even a bit more and require the ordered $\mathrm{K}_{0}$-group of $\mathcal{D}$ to be a subgroup of $\mathbb{Q}$ (with natural order). 
\en

\bn
Unlike the original quasidiagonality question, the version of \ref{finite-ssa-q} yields an obvious infinite counterpart by simply replacing $\mathcal{Q}$ with $\mathcal{Q} \otimes \mathcal{O}_{\infty}$ and `finite' with the minimal necessary condition `not isomorphic to $\mathcal{O}_{2}$':

\bigskip
\noindent 
{${\bf QDQ}_{\mathrm{infinite \; s.s.a.}}$} \enspace {\it Is $\mathcal{D} \otimes \mathcal{Q} \otimes \mathcal{O}_{\infty} \cong \mathcal{Q} \otimes \mathcal{O}_{\infty}$ for every strongly self-absorbing $\mathrm{C}^{*}$-algebra $\mathcal{D}$ not isomorphic to $\mathcal{O}_{2}$?} 

\bigskip
\noindent
Once again this asks for an abstract characterisation of $\mathcal{Q} \otimes \mathcal{O}_{\infty}$ as the final object in the category of all strongly self-absorbing $\mathrm{C}^{*}$-algebras which are not $\mathcal{O}_{2}$ (or equivalently, which have nontrivial $\mathrm{K}$-theory). 

This infinite (or  rather, general) version of the strongly self-absorbing quasidiagonality question runs completely parallel with its finite antagonist, and may be taken as first evidence that the original quasidiagonality question is just the finite incarnation of a much more general type of embedding problem. 
\en

\bn
We have now used a tool from classification---Elliott's intertwining argument---to rephrase the quasidiagonality question as an isomorphism problem, which makes sense both in a finite and an infinite context. Going only one step further, we see that classification not only predicts, but in fact provides, a positive answer to {${\bf QDQ}_{\mathrm{infinite \; s.s.a.}}$}: The secret extra ingredient is to assume that $\mathcal{D}$ satisfies the UCT. Under this hypothesis, it was observed in \cite{TW:TAMS} that the $\mathrm{K}$-theory of $\mathcal{D}$ has to agree with that of one of the known strongly self-absorbing examples, and then it follows from Kirchberg--Phillips classification that  $\mathcal{D}$ is indeed absorbed by $\mathcal{Q} \otimes \mathcal{O}_{\infty}$. We therefore have:

\bigskip
\begin{ntheorem}
If $\mathcal{D} \neq \mathcal{O}_{2}$ is a strongly self-absorbing $\mathrm{C}^{*}$-algebra which satisfies the UCT, then $\mathcal{D} \otimes \mathcal{Q} \otimes \mathcal{O}_{\infty} \cong \mathcal{Q} \otimes \mathcal{O}_{\infty}$.

In other words, $\mathcal{Q} \otimes \mathcal{O}_{\infty}$ is the unique final object in the category of strongly self-absorbing $\mathrm{C}^{*}$-algebras which have nontrivial $\mathrm{K}$-theory and satisfy the UCT.
\end{ntheorem}

\bigskip
\noindent
With this observation at hand, I found it harder and harder to imagine {${\bf QDQ}_{\mathrm{finite \; s.s.a.}}$} fails when also assuming the UCT. Now we know this perception was indeed correct (cf.\ \ref{ssa-cor} below), even in a generality going far beyond the strongly self-absorbing context (see \ref{qd-cor}). Here I took the strongly self-absorbing perspective mostly for a cleaner picture of a simpler situation---but with the benefit of hindsight, the theorem above provided just the necessary impetus to combine the quasidiagonality question with the UCT problem. 
\en

\section{The main result: structure and classification}
\label{main-section}

\bn
\label{main-thm}
\begin{ntheorem}{\rm \cite[Theorem~A]{TWW} }
Let $A$ be a separable nuclear $\mathrm{C}^{*}$-algebra which satisfies the UCT. Then every faithful trace on $A$ is quasidiagonal.
\end{ntheorem}

\bigskip
Short after the distribution of \cite{TWW}, Gabe observed in \cite{Gab:qd-exact} that essentially the same argument works when weakening the nuclearity hypotheses to $A$ being exact and the trace being amenable. Before outlining the proof of the theorem above let us list some consequences, mostly for the structure and classification of simple $\mathrm{C}^{*}$-algebras, but also for Rosenberg's conjecture.
\en

\bn
\label{qd-cor}
\begin{ncor}{\rm \cite[Corollary~B]{TWW} }
Every separable nuclear $\mathrm{C}^{*}$-algebra which satisfies the UCT and has a faithful trace is quasidiagonal. In particular, the quasidiagonality question has a positive answer for simple unital $\mathrm{C}^{*}$-algebras satisfying the UCT.
\end{ncor}
\en

\bn
\label{rosenberg-cor}
In the appendix of \cite{Hadwin:JOT}, Rosenberg observed that for a discrete group $G$, if the reduced group $\mathrm{C}^{*}$-algebra $\mathrm{C}^{*}_{\mathrm{r}}(G)$ is quasidiagonal then $G$ is amenable. The converse was Rosenberg's conjecture, open since the 1980's. Our result in conjunction with \cite{Tu:KT} (which verifies the UCT assumption) confirms the conjecture (the canonical trace on $\mathrm{C}^{*}_{\mathrm{r}}(G)$ is well-known to be faithful). Together with Rosenberg's earlier result this yields a new characterisation of amenability for discrete groups. Note that at first glance our result seems to only cover countable discrecte groups (Theorem~\ref{main-thm} deals with separable $\mathrm{C}^{*}$-algebras), but the general case follows since both quasidiagonality and amenability are \emph{local} conditions.
\bigskip

\begin{ncor} {\rm \cite[Corollary~C]{TWW} }
For a discrete amenable group $G$, its reduced group $\mathrm{C}^{*}$-algebra $\mathrm{C}^{*}_{\mathrm{r}}(G)$ is quasidiagonal.
\end{ncor}
\en

\bn
\label{class-cor}
Elliott, Gong, Lin and Niu have very recently (see \cite{EGLN:arXiv}, which heavily uses \cite{GLN:arXiv}) obtained a spectacular classification result for unital simple nuclear $\mathrm{C}^{*}$-algebras---the crucial additional assumptions being finite decomposition rank and the UCT. They also show that finite decomposition rank may be weakened to finite nuclear dimension, provided all traces are quasidiagonal. Our Theorem~\ref{main-thm} now shows that this last hypothesis is in fact redundant. This is important for applications, since finite nuclear dimension is notoriously easier to verify than finite decomposition rank, but it is also very satisfactory from a conceptual point of view, since for once it allows to state the purely infinite and the stably finite incarnations of classification in the same framework---and it also shows that quasidiagonality of traces precisely marks the dividing line between nuclear dimension and decomposition rank (at least in the simple UCT case), thus answering  \cite[Question~9.1]{WinZac:dimnuc} in this context.

\bigskip

\begin{ncor} {\rm \cite[Corollary~D]{TWW} }
The class of all separable, unital, simple, infinite dimensional $\mathrm{C}^{*}$-algebras with finite nuclear dimension and which satisfy the UCT is classified by the Elliott invariant.
\end{ncor}
\en

\bn
\label{monotracial-cor}
It is worth highlighting the special case when there is at most one trace. For once, the statement becomes particularly clean then, partly because the classifying invariant reduces to just ordered $\mathrm{K}$-theory in this situation, and moreover the proof only relies on work that has already been published (apart from \cite{TWW}). The traceless case has been known for a long time---it is the by now classical Kirchberg--Phillips classification of purely infinite $\mathrm{C}^{*}$-algebras. The equivalence of conditions (i), (ii) and (iii) below in the tracial case is the culmination of \cite{R:IJM,W:Invent2,MS:Acta,SWW:Invent} and does not require the UCT; this only comes in to make the connection with (i'). 

\bigskip

\begin{ncor} {\rm \cite[Corollaries~E and 6.4]{TWW} }
The full Toms--Winter conjecture holds for $\mathrm{C}^{*}$-algebras with at most one trace and which satisfy the UCT. 

That is, for a separable, unital, simple, infinite dimensional, nuclear $\mathrm{C}^{*}$-algebra $A$ with at most one trace and with the UCT, the following are equivalent:
\begin{enumerate}[\hspace{18pt}]
\item[{\rm (i)}] $A$ has finite nuclear dimension.
\item[{\rm (ii)}] $A$ is $\mathcal{Z}$-stable.
\item[{\rm (iii)}] $A$ has strict comparison of positive elements.
\end{enumerate}
If $A$ is stably finite, then {\rm (i)} may be replaced by 
\begin{enumerate}[\hspace{18pt}]
\item[{\rm (i')}]  $A$ has finite decomposition rank.
\end{enumerate}

\noindent
Moreover, this class is classified up to $\mathcal{Z}$-stability by ordered $\mathrm{K}$-theory.
\end{ncor}
\en

\bn
\label{ssa-cor}
Since strongly self-absorbing $\mathrm{C}^{*}$-algebras are $\mathcal{Z}$-stable by \cite{W:JNCG} and have at most one trace, we now know that the chart of \ref{ssa-examples} is indeed complete within the  UCT class.

\bigskip

\begin{ncor}
The strongly self-absorbing $\mathrm{C}^{*}$-algebras satisfying the UCT are precisely the known ones.
\end{ncor}
\en

\section{A sketch of a proof}
\label{outline-section}

\bn
In this outline of the proof of Theorem~\ref{main-thm} I freely assume $A$ to be unital, since one can easily reduce to this case. The very rough idea of the argument is it to produce two complementary cones over $A$ and `connect' them along the interval in order to construct an almost multiplicative map from $\mathcal{C}([0,1]) \otimes A$ to  $M_{2}(\mathcal{Q}_{\omega})$. 
\en

 \bn
Let us begin by producing two cones over $A$ in $\mathcal{Q}_{\omega}$ such that at least their scalar parts are compatible. In order to conjure up a single cone over $A$ inside $\mathcal{Q}_{\omega}$ one might try to employ Voiculescu's theorem \cite{Voi:Duke} on homotopy invariance of quasidiagonality, which will immediately yield an embedding of the cone over $A$ into $\mathcal{Q}_{\omega}$. However, this method will typically give an embedding which is small in trace (not surprisingly, since Voiculescu's result works in complete generality, even when there are no traces around at all). For us this means that we won't be able to repeat the step in order to find the complementary second cone. Instead, we will need a more refined way of implementing quasidiagonality of cones. We will do this by carefully controlling tracial information for the embedding $\mathcal{C}_{0}((0,1],A) \hookrightarrow \mathcal{Q}_{\omega}$. Roughly speaking, we want the canonical trace on $\mathcal{Q}_{\omega}$ to be compatible with a prescribed trace on $A$ and with Lebesgue measure on the interval. This was essentially laid out in \cite{SWW:Invent} and refined in \cite{TWW}; it heavily relies on Connes' \cite{C:Ann} and also uses Kirchberg and R{\o}rdam's \cite{KR:Crelle}.    

\bigskip

\begin{nlemma}
Let $A$ be a separable, unital, nuclear $\mathrm{C}^{*}$-algebra and let $\tau_{A} \in \mathrm{T}(A)$ be a tracial state.
\begin{enumerate}
\item[{\rm (i)}] There is a $^{*}$-homomorphism 
\[
\Psi: \mathcal{C}_{0}((0,1]) \otimes A \to \mathcal{Q}_{\omega}
\]
such that 
\[
\tau_{\mathcal{Q}_{\omega}} \circ \Psi = \mathrm{ev}_{1} \otimes \tau_{A}.
\]
\item[{\rm (ii)}] There are $^{*}$-homomorphisms
\[
\acute{\Phi}: \mathcal{C}_{0}((0,1]) \otimes A \to \mathcal{Q}_{\omega},
\]
\[
\grave{\Phi}: \mathcal{C}_{0}([0,1)) \otimes A \to \mathcal{Q}_{\omega},
\]
\[
\Theta: \mathcal{C}([0,1]) \to \mathcal{Q}_{\omega}
\]
which are compatible in the sense that 
\[
\acute{\Phi}|_{\mathcal{C}_{0}((0,1]) \otimes 1_{A}}= \Theta|_{\mathcal{C}_{0}((0,1])}
\quad
\mbox{ and }  
\quad
\grave{\Phi}|_{\mathcal{C}_{0}([0,1)) \otimes 1_{A}} = \Theta|_{\mathcal{C}_{0}([0,1))},
\]
and such that
\[
\tau_{\mathcal{Q}_{\omega}} \circ \acute{\Phi} = \tau_{\mathrm{\, Lebesgue}}\otimes \tau_{A}
\quad
\mbox{ and }
\quad
\tau_{\mathcal{Q}_{\omega}} \circ \grave{\Phi} = \tau_{\mathrm{\, Lebesgue}}\otimes \tau_{A}.
\]
We use $\tau_{\mathrm{\, Lebesgue}}$ to denote the traces induced by Lebesgue measure on $\mathcal{C}([0,1])$ and on the two cones $\mathcal{C}_{0}((0,1])$ and $\mathcal{C}_{0}([0,1))$.
\end{enumerate}
\end{nlemma}

\bigskip

\begin{idea}
(i) This is essentially contained in \cite{SWW:Invent}. For simplicity let us assume the trace $\tau_{A}$ is extremal, so that the weak closure of the GNS representation of $A$ is a finite injective factor. We therefore obtain a unital $^{*}$-homomorphism $A \to \mathcal{R} \subset \mathcal{R}^{\omega}$ which picks up the trace $\tau_{A}$ when composed with the canonical trace on $\mathcal{R}^{\omega}$. By Kaplansky's density theorem $\mathcal{Q}_{\omega}$ surjects onto $\mathcal{R}^{\omega}$, when dividing out the trace kernel ideal $\{x \in \mathcal{Q}_{\omega} \mid \tau_{\mathcal{Q}_{\omega}}(x^{*}x) = 0\} \lhd \mathcal{Q}_{\omega}$. By the Choi--Effros lifting theorem, there is a c.p.c.\ lift from $A$ to $\mathcal{Q}_{\omega}$. The `curvature' of this map (the defect of it being multiplicative) then lies in the trace kernel ideal of $\mathcal{Q}_{\omega}$, and one can use a quasicentral approximate unit in conjunction with a reindexing argument to replace it by a c.p.c.\ order zero lift. (Alternatively, one can use Kirchberg's $\varepsilon$-test from \cite{KR:Crelle} in place of  reindexing.) This order zero map corresponds to a $^{*}$-homomorphism $\Psi$ defined on the cone over $A$ which will have the right properties.

\bigskip

(ii) Find a $^{*}$-homomorphism $\lambda: \mathcal{C}_{0}((0,1]) \to \mathcal{Q}$ such that $\tau_{\mathcal{Q}} \circ \lambda = \tau_{\mathrm{\, Lebesgue}}$. Next, find a unital copy of $\mathcal{Q}$ in $\mathcal{Q}_{\omega} \cap \Psi(\mathcal{C}_{0}((0,1]) \otimes A)'$. We compose this inclusion with $\lambda$ and tensor with $\Psi$ to obtain a $^{*}$-homomorphism
\[
\widetilde{\Psi}: \mathcal{C}_{0}((0,1]) \otimes \mathcal{C}_{0}((0,1]) \otimes A \to \mathcal{Q}_{\omega}.
\]
Since $\mathcal{C}_{0}((0,1])$ is the universal $\mathrm{C}^{*}$-algebra generated by a positive contraction, the assignment $\mathrm{id}_{(0,1]} \otimes a \mapsto \mathrm{id}_{(0,1]} \otimes \mathrm{id}_{(0,1]} \otimes a$ induces a $^{*}$-homomorphism; we define $\acute{\Phi}$ to be the composition with $\widetilde{\Psi}$.

Next observe that the two cones in $\mathcal{Q}_{\omega}$ generated by the elements $\acute{\Phi}(\mathrm{id}_{(0,1]} \otimes 1_{A})$ and $1_{\mathcal{Q}_{\omega}} - \acute{\Phi}(\mathrm{id}_{(0,1]} \otimes 1_{A})$ carry the same Cuntz semigroup information (which is determined by Lebesgue measure on the interval), and are therefore unitarily equivalent by \cite{CE:IMRN} (and reindexing), i.e., $1_{\mathcal{Q}_{\omega}} - \acute{\Phi}(\mathrm{id}_{(0,1]} \otimes 1_{A}) = w^{*} \, \acute{\Phi}(\mathrm{id}_{(0,1]} \otimes 1_{A}) \, w$ for some unitary $w \in \mathcal{Q}_{\omega}$. Define $\grave{\Phi}$ to be the resulting conjugate of $\acute{\Phi}$, so that $\grave{\Phi}((1 - \mathrm{id}_{[0,1)}) \otimes a) = w^{*} \, \acute{\Phi}(\mathrm{id}_{(0,1]} \otimes a) w$, $a \in A$. The map $\Theta$ is then fixed by these data since $\grave{\Phi}((1 - \mathrm{id}_{[0,1)}) \otimes 1_{A}) + \acute{\Phi}(\mathrm{id}_{(0,1]} \otimes 1_{A})= 1_{\mathcal{Q}_{\omega}}$.
\end{idea}
\en

\bn
\label{matrix-multiplicative}
Now we have produced two cones over $A$ inside $\mathcal{Q}_{\omega}$; the scalar parts of these fit nicely together, but the $A$-valued components might be in general position. The task is to join them in order to find a c.p.c.\ map from $\mathcal{C}([0,1]) \otimes A$ to (matrices over) $\mathcal{Q}_{\omega}$ which is either exactly or at least approximately multiplicative. We wish to establish this  connection by comparing the two restrictions to the suspension over $A$,
\[
\acute{\Lambda}:= \acute{\Phi}|_{\mathcal{C}_{0}((0,1)) \otimes A}: \mathcal{C}_{0}((0,1)) \otimes A \to \mathcal{Q}_{\omega}
\]
and
\[
\grave{\Lambda}:= \grave{\Phi}|_{\mathcal{C}_{0}((0,1)) \otimes A}: \mathcal{C}_{0}((0,1)) \otimes A \to \mathcal{Q}_{\omega}.
\] 
Here's what \emph{would} make this work. It's not quite going to, but it is a blueprint of the actual proof, and it isolates the necessary ingredients. 

\bigskip

\begin{nlemma}
In the setting above, suppose there is a unitary $u \in \mathcal{Q}_{\omega}$ such that
\begin{equation}
\label{Lambda-ue}
\grave{\Lambda}(\, .\,) = u \acute{\Lambda}(\, .\,) \, u^{*}.
\end{equation}
Then, there is a $^{*}$-homomorphism 
\[
\overline{\Phi}: A \to M_{2}(\mathcal{Q}_{\omega})
\]
such that 
\[
(\mathrm{tr}_{M_{2}} \otimes \tau_{\mathcal{Q}_{\omega}}) \circ \overline{\Phi} = \frac{1}{2} \cdot  \tau_{A}.
\]
In particular, the unital $^{*}$-homomorphism 
\[
\widetilde{\Phi}: A \stackrel{\overline{\Phi}}{\longrightarrow} \overline{\Phi}(1_{A}) \, M_{2}(\mathcal{Q}_{\omega}) \, \overline{\Phi}(1_{A}) \cong \mathcal{Q}_{\omega}
\]
shows that the trace $\tau_{A}$ is quasidiagonal.
\end{nlemma}

\bigskip

\begin{nproof}
We write down the map $\overline{\Phi}$ explicitly. Define a partition of unity of piecewise linear positive continuous functions $h_{0},h_{1/2},h_{1}$ on the interval so that $h_{0}$ equals $1$ at $0$, and is $0$ on $[1/4,1]$; $h_{1}$ is just $h_{0}$ reflected at $1/2$. Consider a continuous rotation
$R \in M_{2}(\mathcal{C}([0,1]))$ with 
\[
R|_{[0,1/4]} \equiv \left(\begin{array}{cc} 1 & 0 \\ 0 & 1 \end{array}\right) \quad \mbox{ and } \quad
R|_{[3/4,1]} \equiv \left(\begin{array}{cc} 0 & 1 \\ 1 & 0 \end{array}\right).
\]
Let $\Theta^{(2)} : M_{2}(\mathcal{C}([0,1])) \to M_{2}(\mathcal{Q}_{\omega})$ denote the amplification of $\Theta$ to $2\times2$ matrices.

We may then define a c.p.\ map by setting
\[
\begin{array}{rcl}
\overline{\Phi}(a) & := &  \left( \begin{array}{ccc} \grave{\Phi}(h_{0} \otimes a) & &0 \\  0 && 0 \end{array} \right) \\
\\
&& \quad + \left( \begin{array}{ccc} 1 && 0 \\  0 && u^{*} \end{array} \right) \Theta^{(2)}(R) \left( \begin{array}{ccc}  \grave{\Lambda}(h_{1/2} \otimes a) && 0 \\ 0 && 0\end{array} \right) \Theta^{(2)}(R^{*}) \left( \begin{array}{ccc} 1 && 0 \\  0 && u \end{array} \right) \\
\\
&& \quad \quad + \left( \begin{array}{ccc} 0 & & 0 \\ 0 & &  \acute{\Phi}(h_{1} \otimes a) \end{array} \right)
\end{array}
\]
 for $a \in A$; it is not hard to check that $\overline{\Phi}$ is in fact multiplicative and picks up half of the trace $\tau_{A}$ as claimed in the lemma.
 
For the last statement note that $\overline{\Phi}$ is unital when regarded as a $^{*}$-homomorphism to the hereditary $\mathrm{C}^{*}$-subalgebra generated by its image,  $\overline{\Phi}(1_{A}) \, M_{2}(\mathcal{Q}_{\omega}) \, \overline{\Phi}(1_{A})$, which is isomorphic to $\mathcal{Q}_{\omega}$ since $\mathcal{Q}$ is self-similar; cf.\ Remark~\ref{unital-qd}.  Under this identification the traces $2\cdot (\mathrm{tr}_{M_{2}} \otimes \tau_{\mathcal{Q}_{\omega}})|_{\overline{\Phi}(1_{A}) \, M_{2}(\mathcal{Q}_{\omega}) \, \overline{\Phi}(1_{A})}$ and $\tau_{\mathcal{Q}_{\omega}}$ agree since $\mathcal{Q}_{\omega}$  is monotracial by \cite[Lemma~4.7]{MS:AJM}, so that $\tau_{A} = \tau_{\mathcal{Q}_{\omega}} \circ \widetilde{\Phi}$ is quasidiagonal by Proposition~\ref{qd-traces-prop}.
\end{nproof}
\en

\bn
\label{matrix-multiplicative-remark}
\begin{nremarks}
(i) If one only had an approximate version of (\ref{Lambda-ue}) the same argument would yield an approximately multiplicative c.p.c.\ map $\overline{\Phi}$; after reindexing this would still be good enough to prove quasidiagonality.

(ii) It is natural to ask whether the use of $2 \times 2$ matrices is really essential here. One could certainly hide the matrices by rotating and compressing everything into the upper left corner---but that's a red herring since one cannot necessarily force the resulting map to be unital. The reason is that the method above allows only limited control over $\mathrm{K}$-theory, and one cannot guarantee that $\overline{\Phi}(1_{A})$ is Murray--von Neumann equivalent to $e_{11} \otimes 1_{\mathcal{Q}_{\omega}}$ (of course the two agree tracially, but that's not enough in ultrapowers, even of UHF algebras).    
\end{nremarks}
\en

\bn
In general, unitary equivalence of the two suspensions as in (\ref{Lambda-ue}) seems too much to ask for---and the same goes for approximate versions. On the other hand, it's not completely outrageous either; for example, it is not too hard to see that when $A$ happens to be strongly self-absorbing then the converse of Lemma~\ref{matrix-multiplicative} holds, i.e., unitary equivalence of the two suspensions is implied by quasidiagonality. More can be said using \cite{BBSTWW:arXiv}, but whether this kind of unitary equivalence is a necessary condition for quasidiagonality in complete generality is not clear, and we don't have means to check it directly. The way around this is the stable uniqueness machinery as introduced by Lin in \cite{Lin:stableJOT}, then refined by Dadarlat--Eilers in \cite{DE:PLMS} and since often used and further refined by Elliott, Gong, Lin, Niu, and others.  
\en

\bn
\label{multiplicative-along-interval}
Let us revisit Lemma~\ref{matrix-multiplicative} and replace the critical hypothesis (\ref{Lambda-ue}) by a weaker one (still not quite weak enough for us to confirm it in sufficient generality, but now almost within reach):
\[
\exists \, n \in \mathbb{N}, \, u,v \in \mathcal{U}(M_{n+1}(\mathcal{Q}_{\omega}) )\mbox{ such that } 
\]
\begin{equation}
\label{Lambda-Lambda-ue}
\grave{\Lambda} \oplus \grave{\Lambda}^{\oplus n} = u \Big(\acute{\Lambda} \oplus \grave{\Lambda}^{\oplus n} \Big) u^{*} \quad \mbox{ and } \quad \grave{\Lambda} \oplus \acute{\Lambda}^{\oplus n} = v \Big(\acute{\Lambda} \oplus \acute{\Lambda}^{\oplus n} \Big) v^{*}.
\end{equation}
Then we can chop the interval which sits via $\Theta$ inside $\mathcal{Q}_{\omega} \cong M_{2n} \otimes \mathcal{Q}_{\omega}$ into small pieces and apply the idea of the proof of Lemma~\ref{matrix-multiplicative} $2n$ times along the interval (we have to switch from $u$ to $v$ halfway, which is why we have to use $2n$ intervals, not just $n$); diagrammatically we end up with the following picture; cf.\ \cite[Figure~1]{TWW}:

\begin{picture}(800,180)(40,40)
\setlength{\unitlength}{0.16mm}
\put(82.5,195){\tiny{$\grave\Lambda$}}  
\put(82.5,223){\tiny{$\grave\Lambda$}}  
\put(82.5,251){\tiny{$\grave\Lambda$}}  
\put(82.5,279){\tiny{$\grave\Lambda$}}  
\put(82.5,335){\tiny{$\grave\Lambda$}}  
\put(82.5,363){\tiny{$\grave\Lambda$}}  
\put(82.5,391){\tiny{$\grave\Lambda$}}  
\put(192.5,195){\tiny{$\acute\Lambda$}}  
\put(192.5,223){\tiny{$\grave\Lambda$}}  
\put(192.5,251){\tiny{$\grave\Lambda$}}  
\put(192.5,279){\tiny{$\grave\Lambda$}}  
\put(192.5,335){\tiny{$\grave\Lambda$}}  
\put(192.5,363){\tiny{$\grave\Lambda$}}  
\put(192.5,391){\tiny{$\grave\Lambda$}}  
\put(302.5,195){\tiny{$\acute\Lambda$}}  
\put(302.5,223){\tiny{$\acute\Lambda$}}  
\put(302.5,251){\tiny{$\grave\Lambda$}}  
\put(302.5,279){\tiny{$\grave\Lambda$}}  
\put(302.5,335){\tiny{$\grave\Lambda$}}  
\put(302.5,363){\tiny{$\grave\Lambda$}}  
\put(302.5,391){\tiny{$\grave\Lambda$}}  
\put(412.5,195){\tiny{$\acute\Lambda$}}  
\put(412.5,223){\tiny{$\acute\Lambda$}}  
\put(412.5,251){\tiny{$\acute\Lambda$}}  
\put(412.5,279){\tiny{$\grave\Lambda$}}  
\put(412.5,335){\tiny{$\grave\Lambda$}}  
\put(412.5,363){\tiny{$\grave\Lambda$}}  
\put(412.5,391){\tiny{$\grave\Lambda$}}  
\put(542.5,195){\tiny{$\acute\Lambda$}}  
\put(542.5,223){\tiny{$\acute\Lambda$}}  
\put(542.5,251){\tiny{$\acute\Lambda$}}  
\put(542.5,279){\tiny{$\acute\Lambda$}}  
\put(542.5,335){\tiny{$\acute\Lambda$}}  
\put(542.5,363){\tiny{$\grave\Lambda$}}  
\put(542.5,391){\tiny{$\grave\Lambda$}}  
\put(652.5,195){\tiny{$\acute\Lambda$}}  
\put(652.5,223){\tiny{$\acute\Lambda$}}  
\put(652.5,251){\tiny{$\acute\Lambda$}}  
\put(652.5,279){\tiny{$\acute\Lambda$}}  
\put(652.5,335){\tiny{$\acute\Lambda$}}  
\put(652.5,363){\tiny{$\acute\Lambda$}}  
\put(652.5,391){\tiny{$\grave\Lambda$}}  
\put(762.5,195){\tiny{$\acute\Lambda$}}  
\put(762.5,223){\tiny{$\acute\Lambda$}}  
\put(762.5,251){\tiny{$\acute\Lambda$}}  
\put(762.5,279){\tiny{$\acute\Lambda$}}  
\put(762.5,335){\tiny{$\acute\Lambda$}}  
\put(762.5,363){\tiny{$\acute\Lambda$}}  
\put(762.5,391){\tiny{$\acute\Lambda$}}  
\put(85.5,307){\tiny{$\vdots$}}  
\put(195.5,307){\tiny{$\vdots$}}  
\put(305.5,307){\tiny{$\vdots$}}  
\put(415.5,307){\tiny{$\vdots$}}  
\put(545.5,307){\tiny{$\vdots$}}  
\put(655.5,307){\tiny{$\vdots$}}  
\put(765.5,307){\tiny{$\vdots$}}  
\put(135.5,223){\tiny{$\sim_{u}$}}  
\put(245.5,251){\tiny{$\sim_{u}$}}  
\put(355.5,279){\tiny{$\sim_{u}$}}  
\put(465.5,307){\tiny{$\dots$}}  
\put(595.5,335){\tiny{$\sim_{v}$}}  
\put(705.5,363){\tiny{$\sim_{v}$}}  
\end{picture}
This will produce a $^{*}$-homomorphism 
\[
\overline{\Phi}: A \to M_{2} \otimes M_{2n} \otimes \mathcal{Q}_{\omega}
\]
in a similar way as in Lemma~\ref{matrix-multiplicative}, which again entails quasidiagonality. 
\en

\bn
Just as in Remark~\ref{matrix-multiplicative-remark}(i), it would be enough to come up with an approximately multiplicative c.p.c.\ map $\overline{\Phi}$, which would follow from an approximate version of (\ref{Lambda-Lambda-ue}). The latter is very close to the conclusion of Theorem~\ref{su-original}, with $\grave{\Lambda}$ and $\acute{\Lambda}$ in place of $\varphi$ and $\psi$, respectively, and also with $\grave{\Lambda}$ and $\acute{\Lambda}$ in place of $\iota$. However, there is a catch: The maps in the diagrammatic chart of \ref{multiplicative-along-interval} are in fact not the original maps $\grave{\Lambda}$ or  $\acute{\Lambda}$; rather, they are restrictions of those maps to small subintervals of $(0,1)$. This makes a difference, since it means that the maps depend on the number of intervals, hence on $n$, which in Theorem~\ref{su-original} in turn depends on the maps---and the whole affair becomes circular! Luckily, there is a backdoor: In the controlled stable uniqueness theorem \ref{su-controlled} the number $n$ does not depend on the actual maps; it only depends (except for $\mathcal{G}$ and $\delta$, of course) on the control function $\Delta$ which is tied to the Lebesgue measure on the interval via the prescribed trace and the map $\Theta$. The price for this additional control is the UCT hypothesis in Theorem~\ref{main-thm}.
\en

\section{Some open problems}
\label{open-section}

\bn
Of course the main problems in the context of this paper are the UCT problem and the  quasidiagonality question in its various versions as discussed in Section~\ref{QDQ-section}. These are expected to be hard; the problems listed below aim to highlight their interplay and to break them up into smaller bits and pieces which will hopefully be easier to attack. 
\en

\bn
\begin{nquestion}
Are there formal implications between the versions of the quasidiagonality question from Section~\ref{QDQ-section}? 
In other words, can we prove any of the implications  [${\bf QDQ}_{\mathrm{infinite \; s.s.a.}}$ holds] $\Longleftrightarrow$ [${\bf QDQ}_{\mathrm{finite \; s.s.a.}}$ holds] $\Longrightarrow$ [${\bf QDQ}_{\mathrm{simple,1}}$ holds] $\Longrightarrow$ [${\bf QDQ}_{\mathrm{simple}}$ holds] $\Longrightarrow$ [${\bf QDQ}$ holds]?
\end{nquestion}
\en

\bn
By Corollary~\ref{qd-cor}, the UCT implies quasidiagonality under suitable conditions, and one can ask under which hypotheses there is a converse. This is also interesting for special cases:
\bigskip

\begin{nquestions}
Does every quasidiagonal strongly self-absorbing $\mathrm{C}^{*}$-algebra satisfy the UCT?  
What about strongly self-absorbing $\mathrm{C}^{*}$-subalgebras of $ \mathcal{Q} \otimes \mathcal{O}_{\infty}$? Or unital, simple, nuclear and monotracial $\mathrm{C}^{*}$-subalgebras of $\mathcal{Z}$? \end{nquestions}
\en

\bn
Kirchberg has reduced the UCT problem to the simple case; even more, the problem is equivalent to the question whether $\mathcal{O}_{2}$ is the only unital Kirchberg algebra with trivial $\mathrm{K}$-theory (see \cite[2.17]{K:Abel}). As discussed in \ref{ssa-initial-final}, for strongly self-absorbing $\mathrm{C}^{*}$-algebras the answer is known. From this point of view the following does not seem likely, but I still think it is worth asking.

\bigskip

\begin{nquestion}
 Can the UCT problem be reduced to the strongly self-absorbing case? 
\end{nquestion}
\en

\bn
It was shown in \cite{TW:TAMS} that the $\mathrm{K}$-theory of a strongly self-absorbing $\mathrm{C}^{*}$-algebra satisfying the UCT has to agree with the $\mathrm{K}$-theory of one of the known strongly self-absorbing examples. However, the proof really only requires the formally weaker K\"unneth Theorem for tensor products (see \cite{RS:DMJ}), and one may ask whether even this can be made redundant, or whether there are at least some restrictions on the possible $\mathrm{K}$-groups. For example, Dadarlat pointed out that for a quasidiagonal strongly self-absorbing $\mathcal{D}$, $\mathrm{K}_{1}(\mathcal{D})$ cannot have an infinite cyclic subgroup (again by the K\"unneth Theorem and since in this case $\mathcal{D} \otimes \mathcal{Q} \cong \mathcal{Q}$). 

\bigskip
 
\begin{nquestions}
If $\mathcal{D}$ is a strongly self-absorbing $\mathrm{C}^{*}$-algebra, does $\mathrm{K}_{1}(\mathcal{D})$ have to be trivial? Does $\mathrm{K}_{*}(\mathcal{D})$ have to be torsion free?
\end{nquestions}
\en

\bn
It is a classical question when a $\mathrm{C}^{*}$-algebra is isomorphic to its opposite. Whenever one expects classification by $\mathrm{K}$-theory data the answer should be positive, and it certainly is for strongly self-absorbing $\mathrm{C}^{*}$-algebras with UCT; see Corollary~\ref{ssa-cor} (note that the opposite of a strongly self-absorbing $\mathrm{C}^{*}$-algebra is again strongly self-absorbing).  

\bigskip

\begin{nquestion}
Is a strongly self-absorbing $\mathcal{D}$ isomorphic to its opposite $\mathcal{D}^{\mathrm{op}}$?
\end{nquestion}
\en

\bn
The following stems essentially from \cite{B:MAMS}.

\bigskip

\begin{nquestions}
For a separable unital $\mathrm{C}^{*}$-algebra $A$, do the quasidiagonal traces form a face? If, in addition, $A$ is quasidiagonal, are all traces quasidiagonal? Do nuclearity of $A$ or amenability of the traces make a difference?
\end{nquestions}

\bigskip
\noindent
Together with a result from \cite{B:MAMS}, \cite[Corollary~6.1]{TWW} yields a positive answer to the second question when also assuming nuclearity and the UCT.
\en

\bn
In both \cite{ORS:GAFA} and \cite{TWW}, quasidiagonality of amenable group $\mathrm{C}^{*}$-algebras is derived abstractly from classification techniques---but at this point there is no way to construct quasidiagonalising finite-dimensional subspaces of $\ell^{2}(G)$ explicitly.    

\bigskip

\begin{nquestion}
Is there a group theoretic / dynamic proof of Rosenberg's conjecture?
\end{nquestion}
\en

\bn
$\mathrm{C}^{*}$-algebras of amenable groups are almost never simple---but they have simple quotients, and one may ask when these are classifiable. There is by now a range of very convincing results along these lines; cf.\ \cite{ELPW, EckGil:irr}. 

\bigskip

\begin{nquestion}
When are simple quotients of amenable group $\mathrm{C}^{*}$-algebras classifiable? When can one at least show $\mathcal{Z}$-stability?
\end{nquestion}  
\en

\bn
In a similar vein, one can look at topological dynamical systems, where free and minimal actions typically yield simple $\mathrm{C}^{*}$-crossed products. These algebras tend to be nuclear provided the groups---or at least their actions---are amenable; cf.\ \cite{TW:GAFA, Anan:amen}. We know from \cite{GioKer:Crelle} that one cannot expect regularity in general, and that conditions on the dimension of the underlying space (or again the action) are essential to guarantee $\mathcal{Z}$-stability or finite nuclear dimension of the crossed product. Recent results of Kerr, however, together with the tiling result of \cite{DowHucZha}, suggest that we might be only a stone's throw away from an answer to the following:

\bigskip

\begin{nquestion}
For free minimal actions of countable discrete amenable groups on Cantor sets, are the crossed product $\mathrm{C}^{*}$-algebras classifiable?
\end{nquestion}

\bigskip
\noindent
The setup is shockingly general: free minimal Cantor actions of amenable groups! So how are we even entitled to ask this? Quasidiagonality of the crossed product is given by our Theorem~\ref{main-thm} in connection with \cite{Tu:KT}, which verifies the UCT. $\mathcal{Z}$-stability seems now within reach with Kerr's techniques on tiling (based on \cite{DowHucZha}) together with Archey and Phillips' large subalgebra approach \cite{ArPhi:large} or, alternatively, using the idea of dynamic dimension and dynamic $\mathcal{Z}$-stability as defined by the author (yet unpublished, but closely related to the notion of Rokhlin dimension from \cite{HirWinZac:Rokhlin-dimension}).  From here only finite nuclear dimension of the crossed product would be missing to arrive at classifiability (by \cite{EGLN:arXiv} via \cite{Win:AJM}; for a slightly more direct approach in the uniquely ergodic case see \cite{SWW:Invent}). In the case when the ergodic measures form a compact space, finite nuclear dimension follows from $\mathcal{Z}$-stability by \cite{BBSTWW:arXiv}.  

\bigskip
\noindent
Here is an even more general---though not necessarily more daring---layout.

\bigskip
\begin{nquestions}
For free minimal actions of countable discrete amenable groups on finite dimensional spaces, are the crossed product $\mathrm{C}^{*}$-algebras classifiable? What about amenable actions of countable discrete groups?
\end{nquestions}

\bigskip
\noindent
The following rigidity question was beautifully answered for Cantor minimal $\mathbb{Z}$-actions in \cite{GPS:orbit} in terms of strong orbit equivalence. In the situation of amenable group actions, it seems much more speculative, and one can only expect a less complete answer. If one is prepared to go beyond the context of amenable group actions, Popa's rigidity theory for von Neumann algebras (cf.\ \cite{Popa:ICM}) is extremely encouraging---but on the $\mathrm{C}^{*}$-algebra side one would have to change the game completely and develop most of the technology from scratch.

\bigskip

\begin{nquestion}
 To what extent are topological dynamical systems determined by their associated $\mathrm{C}^{*}$-algebras?
\end{nquestion}
\en


\end{document}